\documentclass[12pt]{article}
\usepackage{amsmath,amsxtra,amssymb,latexsym, amscd}
\usepackage[mathscr]{eucal}

\begin{document}                        
                        \def\DATE{1/16/'06}
\title{ {\huge\bf Persistence and \\
global attractivity\\
in the model $A_{n+1}=A_nF(A_{n-m})$} }
\def\1{\rule{0cm}{0cm}} \def\qd{\rule{3mm}{3mm}} \def\BB{$\bullet$}
\renewcommand{\arraystretch}{1.25}
\renewcommand{\theequation}{\thesectn.\arabic{equation}}
\def\sce{\setcounter{equation}{0}}  \newcounter{sectn} \newcounter{sbsect}
\def\sect#1{\addtocounter{section}{1}\sce\setcounter{sbsect}{0}%
        \renewcommand{\thesectn}{\thesection}\1\smallskip\\
        {\1\hspace{-2em}\large\bf\thesectn.\qquad #1\smallskip\par}}
\def\subsect#1{\addtocounter{sbsect}{1}\sce%
        \renewcommand{\thesectn}{\thesection:\Alph{sbsect}}\1\smallskip\\
        {\bf\1\hspace{-1.5em}\thesectn.\qquad #1\smallskip\par}}
\def\THM#1#2{\1\smallskip\par\noindent{\bf THEOREM#1:}
                                        \qquad{\sl #2}\smallskip\\}
\def\LEM#1#2{\1\smallskip\par\noindent{\bf LEMMA#1:}\qquad{\sl #2}\smallskip\\}
\def\COR#1#2{\1\smallskip\par\noindent{\bf COROLLARY#1:}
                                        \qquad{\sl #2}\smallskip\\}
\def\proof{\bigskip\noindent {\sc Proof:}\qquad}
\def\REM#1#2{\1\smallskip\par\noindent{\bf REMARK#1 }\qquad #2\smallskip\\}
\def\qed{\hfill$\quad$\qd\medskip\\} \def\ds{\displaystyle}
\def\LB#1{\label{#1}} \def\BE#1#2{\begin{#1} #2 \end{#1}}
\def\EQ#1#2{\BE{equation}{\LB{#1} #2}} \def\ARR#1#2{\BE{array}{{#1} #2}}
\def\DES#1{\BE{description}{#1}} \def\QT#1{\BE{quote}{#1}}
\def\ENUM#1{\BE{enumerate}{#1}} \def\ITM#1{\BE{itemize}{#1}}
 \def\COM#1{\par\noindent{\bf COMMENT:\quad\sl #1}\par\noindent}
\def\mapsfrom{\hbox{$\;{\leftarrow}\kern-.15em{\mapstochar}\:\:$}}
\def\vv{\kern.344em{\rule[.18ex]{.075em}{1.32ex}}\kern-.344em}
\def\RE{\mbox{\rm I\kern-.21em R}} \def\CX{\mbox{\rm \vv C}}
\def\imp{\Rightarrow} \def\emb{\hookrightarrow} \def\wk{\rightharpoonup}
\def\rd{\dot{\1}} \def\d{\cdot} \def\+{\oplus} \def\x{\times}
\def\<{\langle} \def\>{\rangle} \def\o{\circ} \def\at#1{\Bigr|_{#1}}
\def\cd{\partial} \def\grad{\nabla} \def\L{\left} \def\R{\right}
\def\eps{\varepsilon} \def\f{\varphi} \def\om{\omega} \def\Om{\Omega}
\def\gm{\gamma} \def\gep{\gm_\eps} \def\lm{ } \def\lep{\lm_\eps}
\def\dl{\delta} \def\rb{\bar{r}} \def\rt{\tilde{r}} \def\al{\alpha}
\def\et{e^{-(m+1)}} \def\fh{\hat{f}}
\def\bx{\mathbf{x}} \def\by{\mathbf{y}} \def\bS{\mathbf{S}}
\def\H{{\cal H}} \def\U{{\cal U}} \def\D{{\cal D}}\def\bc{\mathbf{c}}
\def\eq{equation} \def\de{differential \eq} \def\pde{partial \de}
\def\sol{solution} \def\pb{problem} \def\bdy{boundary} \def\fn{function}
\def\dde{delay \de} \def\ev{eigenvalue}
\def\R{\mathbf R}
\author{
{\bf Dang Vu Giang}\\
Hanoi Institute of Mathematics\\
18 Hoang Quoc Viet, 10307 Hanoi, Vietnam\\
{\footnotesize          e-mail: $\<$dangvugiang@yahoo.com$\>$}\\
 }
\maketitle

\noindent {\bf Abstract.}  First, we systemize ealier results the uniform persistence for
 discrete model $A_{n+1}=A_nF( A_{n-m})$ of population growth, where $F:(0,\infty)\to(0,\infty)$ 
is continuous and strictly decreasing.   Second, we investigation the effect of
delay $m$  when $F$ is not monotone. We are mainly using $\omega$-limit set of
persistent solution, which is discussed in more general by P.
Walters, 1982.

\bigskip



\eject

\sect{Introduction }
\bigskip
Recently, several authors (Graef et al 1998, 2006, Liz 2007, Tkachenko et al 2006)  study  the model 
\EQ{a1} {A_{n+1}=A_n F(A_{n-m}),\qquad n=0,1,\cdots, }
where $F:(0,\infty)\to(0,\infty)$ is continuous and strictly decreasing.  They assume that there is 
a unique positive equilibrium $\bar x$ of (\ref{a1}).  Clearly, $F(\bar x)=1$. We say that $\bar x$ is globally attractive 
if all solution of (\ref{a1})
convergence to $\bar x$ as $n\to\infty$. Graef et al  (1998) investigate the bobwhite quail population  with
$$F(x)=\alpha+\frac \beta{1+x^r}.$$
They prove that if 
$$r<\frac{2\alpha+\beta+2\sqrt{\alpha^2+\alpha\beta}}{\beta}\frac{3m+4}{2(m+1)^2},$$
then 
$$\bar x=\root r\of{\frac{\alpha+\beta-1}{1-\alpha}}$$
 is global attractive. Liz (2007) improves this result. More exactly, the equilibrium  $\bar x$ is globally
attractive if 
$$r<\frac{\beta}{(\alpha+\beta-1)(1-\alpha)}\frac{3m+4}{2(m+1)^2}.$$
If $m=0$ this result is sharp. 

\bigskip

\sect{The Persistence}
\bigskip
In this section we assume $F:(0,\infty)\to(0,\infty)$ is continuous only.  
A positive solution $\{A_n\}_{n=-m}^\infty$ is called persistent if
$$0<\liminf_{n\to\infty} A_n \leqslant  \limsup_{n\to\infty} A_n<\infty.$$
The  following theorem gives a sufficient condition for persistent (non-extinctive) populations.
\bigskip

\par{\bf Theorem  1.} {\it Assume that 
\EQ{f0}{0<F(x)\leqslant c<\infty }
 for all $x\ge0$ and
\EQ{f1}{\limsup_{x\to\infty}{F(x)}<1,  }
\EQ{f2}{\liminf_{x\to 0+0}{F(x)}>1.  }
Then every solution  $\{A_n\}_{n=-m}^\infty$ of} (1.1) {\it is persistent.}
\bigskip
\par {\sl Proof}: First, we prove that $\{A_n\}_{n=-m}^\infty$  is
bounded from above. Assume, for the sake of a contradiction, that
$\limsup A_n= \infty$. For each integer $n\ge m$, we define
$$k_n:=\max\{\rho: -m\leqslant \rho\leqslant n, A_\rho=\max_{-m\leqslant i\leqslant n} A_i\}.$$
Observe that $k_{-m}\leqslant  k_{-m+1}\leqslant  \cdots \leqslant  k_n \to\infty$ and that
\EQ{ainfty}{\lim_{n\to\infty} A_{k_n}=\infty.}
Let $n_0>0$ such that $k_{n_0}>0$. We have for $n>n_0$,
$$ A_{k_n-1} F(A_{k_n-1-m})=A_{k_n}\ge  A_{k_n-1}$$
and therefore,
$$ F(A_{k_n-1-m})\ge1.$$
By (\ref{f1}), this implies that
\EQ{afinite}{\limsup_{n\to\infty} A_{k_n-1-m}<\infty.  }
On the other hand,
\begin{eqnarray*}
   A_{k_n} &=&A_{k_n-1}F(A_{k_n-1-m})=\cdots \\
 &=&A_{k_n-1-m}F(A_{k_n-1-2m})\cdots F(A_{k_n-1-m})\\
&\leqslant &\ A_{k_n-1-m}c^{m+1}.
  \end{eqnarray*}
Now take $\limsup$ on both sides we have $\limsup_{n\to\infty} A_{k_n}<\infty$ which
  contradicts (\ref{ainfty}). Thus, $\{A_n\}_{n=-m}^\infty$  is bounded from above. Let $K$ be an upper bound of 
$\{A_n\}_{n=-m}^\infty$ .
\par Next, we prove that $\liminf_{n\to\infty} A_n>0$.
Assume, for the sake of a contradiction, that $\liminf A_n=0$. For
each integer $n\ge m$, we define
$$s_n:=\max\{\rho: -m\leqslant \rho\leqslant  n, A_\rho=\min_{-m\leqslant  i\leqslant  n} A_i\}.$$
Clearly, $s_{-m}\leqslant  s_{-m+1}\leqslant  \cdots \leqslant  s_n \to\infty$ and that
\EQ{ais0}{\lim_{n\to\infty} A_{s_n}=0.  }
Let $n_0>0$ such that $s_{n_0}>0$.
We have for $n>n_0$,
$$ A_{s_n}= A_{s_n-1} F(A_{s_n-1-m})\ge A_{s_n} F(A_{s_n-1-m}) 
$$
and therefore, 
$$ F(A_{s_n-1-m})\leqslant 1.$$
By (\ref{f2}), this implies that
$$\liminf_{n\to\infty} A_{s_n-1-m}>0.$$
On the other hand,
\begin{eqnarray*}
A_{s_n}&=&A_{s_n-1} F(A_{s_n-1-m})=\cdots\\
&=& A_{s_n-1-m}F(A_{s_n-1-2m})\cdots F(A_{s_n-1-m}) \\
&\ge& A_{s_n-1-m}[\min_{x\in[0,K]}F(x)]^{m+1}. 
  \end{eqnarray*}
Now take $\liminf$  as $n\to\infty$ on both sides we have $\liminf_{n\to\infty} A_{s_n}>0$ which contradicts (\ref{ais0})
 The proof is complete.
\bigskip

Next, we will use full-limiting sequences to estimate the $\liminf$ and $\limsup  A_n$.  
Let $\{A_n\}_{n=-m}^\infty$  be a persistent solution of (1.1). Then there are two sequences
$\{P_n\}_{n=-\infty}^\infty$
and
$\{Q_n\}_{n=-\infty}^\infty$
called full  limiting sequences
satisfying equation (1.1) for all $n=0,\pm1,\pm2,\cdots$ such that
$$\limsup_{n\to\infty} A_n=P_0,  \qquad\liminf_{n\to\infty} A_n=Q_0,  $$
and
$$Q_0\leqslant  P_n\leqslant  P_0,\quad Q_0\leqslant  Q_n\leqslant  P_0, \quad\text{ for all }\quad n=0,\pm1,\pm2,\cdots.
$$ The existence of these sequences is guaranteed by Giang et al (2005) and Walter (1982).
It is very convenient to put
$$\alpha=\inf_{x>0}F(x)\qquad\text{ and }\qquad c=\sup_{x>0}F(x).$$
We assume that  there is 
a unique positive equilibrium $\bar x$ of (\ref{a1}) and
\EQ{barx}{0<\alpha< 1=F(\bar x)<c<\infty.}
The following theorem will  give the uniform persistence of the population in model (1.1). 

\bigskip

\par{\bf Theorem  2.} {\it Assume that $(2.1)-(2.3)$ and  $(\ref{barx})$ hold. Then for every solution 
$\{A_n\}_{n=-m}^\infty$ of  $(1.1)$, we have}
\[\bar x\alpha^{m+1}<\liminf_{n\to\infty}A_n\leqslant\bar x\leqslant\limsup_{n\to\infty}A_n
<\bar x c^{m+1}.\]

\bigskip

{\sl Proof:} As the above, let $\{P_n\}_{n=-\infty}^\infty$
and $\{Q_n\}_{n=-\infty}^\infty$ be full time solutions of (1.1) with $P_0=\limsup_{n\to\infty} A_n$ and  $Q_0=\liminf_{n\to\infty} A_n$. We have 
$$ Q_0  =  Q_{-1}F(Q_{-1-m}) \ge   Q_0F(Q_{-1-m}),
$$
so $1\ge F(Q_{-1-m})$.  Our assumptions assure that $F(x)\leqslant 1$ iff $x\ge\bar x$. Hence, $Q_{-1-m}\ge\bar x$, and consequently, $P_0\ge\bar x$.  On the other hand,
\begin{eqnarray*}
Q_0&=&Q_{-1}  F(Q_{-1-m})=Q_{-2}  F(Q_{-2-m}) F(Q_{-1-m})=\cdots\\
&=& Q_{-1-m}  F(Q_{-1-2m})F(Q_{-2m})\cdots F(Q_{-1-m})>  \alpha^{m+1}\bar x.
\end{eqnarray*}
Similarly,
$$ P_0  =   P_{-1}F(P_{-1-m}) \leqslant    P_0F(P_{-1-m}),
$$
so $1\leqslant F(P_{-1-m})$.  But our assumptions imply that $F(x)\ge 1$ if and only if $x\leqslant\bar x$.
Therefore, $P_{-1-m}\leqslant \bar x$. This
implies that $Q_0\leqslant\bar x$. On the other hand, 
\begin{eqnarray*}
P_0&=&P_{-1}  F(P_{-1-m})=P_{-2}  F(P_{-2-m}) F(P_{-1-m})=\cdots\\
&=& P_{-1-m}  F(P_{-1-2m})F(P_{-2m})\cdots F(P_{-1-m})< c^{m+1}\bar x.
\end{eqnarray*}
The proof is now complete.

\bigskip

\sect{3/2-Stabilitiy}

\bigskip

We assume further that 
\EQ{lips}{|\ln F(x)|\leqslant L|\ln(x/\bar x)|\qquad\hbox{ for all }\quad x\in (0,\bar x c^{m+1}).}
\bigskip

{\bf Theorem 3. } {\it Assume that $(2.1)-(2.3)$, $(\ref{barx})$ and $(\ref{lips})$   hold. } {\it Suppose further that}
$$(m+\frac32) L<\frac32.$$
{\it Then every  solution $\{A_n\}_{n=-m}^\infty$ of} (1.1) {\it converges to $\bar x$}.

\bigskip

{\sl Proof:} By Theorem 2, every nonoscillated solution converges to $\bar x$. 
Therefore, without loss of generality we assume that $L(m+\frac32)\ge 1$ and $\{A_n\}_{n=-m}^\infty$
 is an oscillated solution 
(if $L$ is small we can enlarger it). This means that there is a sequence $t_n\to\infty$ of integers such that 
$A_{t_n}\leqslant\bar x$, $A_{t_n+1}>\bar x$  and $t_{n+1}-t_n>3m$ for every $n=1,2,\cdots$. Let 
\[\rho_n>\bigl|\ln\frac{A_t}{\bar x}\bigr|\quad
\hbox{ for every }
\quad t\ge t_n-2m.
 \]
Then 
\[\Bigl|\ln\frac{A_{t+1}}{A_t}\Bigr|=|\ln F(A_{t-m})|\leqslant L\bigl|\ln\frac{A_{t-m}}{\bar x}\bigr|\leqslant L\rho_1\]
for all $t\ge t_1-m$. Let $A_{t_*}\leqslant\bar x$, with $t_*\ge t_1$.
It follows that
\[\Bigl|\ln\frac{A_{s}}{\bar x}\Bigr|\leqslant\sum_{t=s}^{t_*}\bigl|\ln\frac{A_{t+1}}{A_t}\bigr|\leqslant L\rho_1
(t_*+1-s)\]
for all $s\in[t_1-m,t_*]$. Furthermore,
\[\Bigl|\ln\frac{A_{t+1}}{A_t}\Bigr|=|\ln F(A_{t-m})|\leqslant L\bigl|\ln\frac{A_{t-m}}{\bar x}\bigr|\leqslant L^2\rho_1(t_*+m+1-t)\]
for all $t\in[t_1,t_*+m]$.
First, we prove that 
\[\bigl|\ln\frac{A_t}{\bar x}\bigr|\leqslant\rho_1\bigl(L(m+\frac32)-\frac12\bigr)
\quad\hbox{ for  all }\quad t>t_1+m.
\]
If this were not so, let
\[T=\min\Bigl\{t>t_1+m:\quad A_t>\bar x
\quad\bigl|\ln\frac{A_t}{\bar x}\bigr|> \rho_1\Bigl(L(m+\frac32)-\frac12\Bigr)
\Bigr\}\]
If $A_{t_*}:=A_{T-(m+1)}\leqslant\bar x$ then 
\begin{eqnarray*}
|\rho_1\bigl(L(m+\frac32)-\frac12\bigr)|&<&\bigl|\ln\frac{A_T}{\bar x}\bigr| <\sum_{t=t_*}^{t_1+m}
\Bigl|\ln\frac{A_{t+1}}{A_t}\Bigr|\\
&\leqslant&\sum_{t=t_*}^{t_*+m-[\frac1L]}L\rho_1+\sum_{t=t_*+m-[\frac1L]+1}^{t_*+m}L^2\rho_1(t_*+m+1-t)\\
&\leqslant&L\rho_1\bigl(m+1-[\frac1L\bigr])+\frac12\rho_1L^2[\frac1L]([\frac1L]+1)\\
&\leqslant&\rho_1\bigl(L(m+\frac32)-\frac12\bigr).\\
\end{eqnarray*}
This is a contradiction, so we have $A_{T-(m+1)}>\bar x$. Hence, $F(A_{T-(m+1)})<1$ and consequently, $A_{T-1}>A_T$.
By the minimality of $T$, $T=t_1+m+1$, and we have 
\begin{eqnarray*}
|\rho_1\bigl(L(m+\frac32)-\frac12\bigr)|&<&\bigl|\ln\frac{A_T}{\bar x}\bigr| <\sum_{t=t_1}^{t_1+m}
\Bigl|\ln\frac{A_{t+1}}{A_t}\Bigr|\\
&\leqslant&\sum_{t=t_1}^{t_1+m-[\frac1L]}L\rho_1+\sum_{t=t_1+m-[\frac1L]+1}^{t_1+m}L^2\rho_1(t_1+m+1-t)\\
&\leqslant&L\rho_1\bigl(m+1-[\frac1L\bigr])+\frac12\rho_1L^2[\frac1L]([\frac1L]+1)\\
&\leqslant&\rho_1\bigl(L(m+\frac32)-\frac12\bigr).\\
\end{eqnarray*}
This is a contradiction, so we have
\[\bigl|\ln\frac{A_t}{\bar x}\bigr|\leqslant\rho_1\bigl(L(m+\frac32)-\frac12\bigr)
\quad\hbox{ for  all }\quad t>t_1+m.
\]
This result ensures us to choice 
\[\rho_2=\rho_1\bigl(L(m+\frac32)-\frac12\bigr).
\]
Repeat the above argument  (with $t_1$ and $\rho_1$ replaced by $t_2$ and $\rho_2$) 
 we have
\[\bigl|\ln\frac{A_t}{\bar x}\bigr|\leqslant\rho_2\bigl(L(m+\frac32)-\frac12\bigr)
\quad\hbox{ for  all }\quad t>t_2+m.
\]
By induction on $n$, we have
\[\bigl|\ln\frac{A_t}{\bar x}\bigr|\leqslant\rho_1\bigl(L(m+\frac32)-\frac12\bigr)^n
\quad\hbox{ for  all }\quad t>t_n+m.
\]
Using the assumption $\bigl(L(m+\frac32)-\frac12\bigr)<1$, we complete the proof.
\bigskip

\bigskip
\bigskip

\sect{ Application}

\bigskip

Now  consider  again the bobwhite quail population. Here,
$$F(x)=\alpha+\frac \beta{1+x^r}.$$
We can compute 
\[L=\frac{\beta r}{2\alpha+\beta+2\sqrt{\alpha^2+\alpha\beta}}.\]
Hence, if 
$$m+\frac32<\frac{2\alpha+\beta+2\sqrt{\alpha^2+\alpha\beta}}{\beta r}\cdot\frac32,$$
then $\bar x$ is globally attractive.

\bigskip

Next we consider the equation given by Pielou (1974)
\[A_{n+1}=\frac{\beta A_n}{1+\lambda A_{n-m}}.\]
Here 
\[F(x)=\frac\beta{1+\lambda x},\]
and 
\[\bar x=\frac{\beta-1}{\lambda}.\]
We can compute 
\[L=\frac1{\lambda\bar x}=\frac1{\beta-1}.\]
Hence, if 
\[m+\frac32<\frac{3(\beta-1)}2,\]
then $\bar x$ is globally attractive.

\bigskip
\bigskip

\centerline{\huge\bf References}

\bigskip
\par\noindent Giang, D.V., Huong, D.C., 2005. {\it Extinction, Persistence and Global stability in
models of population growth,}  {\sl J. Math. Anal. Appl.} {\bf 308}, pp 195-207.

\bigskip

 \par\noindent  Graef, J.R. , Qian C., Spikes, P.W. ,1998. {\it Stability in a population
model, } {\sl Appl. Math. Comput. } {\bf 89}, pp 119-132.

\bigskip

 \par\noindent Graef, J.R. , Qian C., 2006.  {\it Global attractivity in a nonlinear  difference equation and application, } {\sl Dynam. Systems Appl.} {\bf 15}, pp 89-96.

\bigskip

 \par\noindent  Liz, E.,  2007.  {\it A sharp global stability result for  a discrete  population model,} {\sl J. Math. Anal. Appl.} {\bf 330}, pp 740-743.
\bigskip

 \par\noindent Liz, E.,  Tkachenko, V.,  Trofimchuk, S., 2006. {\it Global stability in discrete  population models 
with delayed-sensity dependence, } {\sl Math. Biosci.} {\bf 199}, pp 26-37.
\bigskip

\par\noindent Pielou, E.C., 1974. "Population and Community Ecology", New York.
\bigskip

 \par\noindent Tkachenko, V.,  Trofimchuk, S., 2006. {\it A global attractivity criterion for nonlinear non-au- tonomous difference equation,} {\sl J. Math. Anal. Appl.} {\bf 322}, pp 901-912.
\bigskip

 \par\noindent Walters, P.,  1982. "An introduction to ergodic theory", Springer-Verlag.

\end{document}